# A Green Perspective on Capacitated Time-dependent Vehicle Routing Problem with Time Windows

## Iman Kazemian

Department of Industrial Engineering
College of Engineering, University of Tehran
North Kargar, Tehran, PO Box: 4563-11155, Iran
E-mail: i.kazemian@ut.ac.ir

## Samin Aref

Department of Computer Science, University of Auckland
Auckland, Private Bag 92019, New Zealand
Computer Science Practice Pathway Group, Unitec Institute of
Technology, Auckland, Private Bag 92025, New Zealand
E-mail: sare618@aucklanduni.ac.nz

**Abstract:** This study presents a novel approach to the vehicle routing problem by focusing on greenhouse gas emissions and fuel consumption aiming to mitigate adverse environmental effects of transportation. A time-dependent model with time windows is developed to incorporate speed and schedule in transportation. The model considers speed limits for different times of the day in a realistic delivery context. Due to the complexity of solving the model, a simulated annealing algorithm is proposed to find solutions with high quality in a timely manner. Our method can be used in practice to lower fuel consumption and greenhouse gas emissions while total route cost is also controlled to some extent. The capability of method is depicted by numerical examples productively solved within 3.5% to the exact optimal for small and mid-sized problems. Moreover, comparatively appropriate solutions are obtained for large problems in averagely one tenth of the exact method restricted computation time.

**Keywords:** Vehicle Routing; Green; Time-dependent; Time windows; Greenhouse gas emissions; Fuel consumption; Simulated Annealing; Meta-heuristic; Mathematical modeling; Optimization

**Biographical notes:** Iman Kazemian holds a B.S. degree in Industrial Engineering from Isfahan University of Technology as well as an M.S. degree from Sharif University of Technology. He is an Industrial Engineering doctoral candidate at University of Tehran. His research interests include Mathematical Programming, Operations Management, and Logistics.

Samin Aref holds an M.S. degree in Industrial Engineering from Sharif University of Technology. He is a Computer Science doctoral candidate at University of Auckland and a lecturer at Unitec Institute of Technology. His fields of interest include Complex Social and Economic Systems, Operations Research, and Structural Analysis of Social Networks.

Authors have presented at New Zealand Mathematics Colloquium, International School and Conference on Network Science, Operations Research Society of New Zealand Conference, and Conference of Iranian Operations Research Society. They have published in *Global Journal of Flexible Systems Management, International Journal of Occupational Safety and Ergonomics, European Journal of Operational Research*, and *Journal of Revenue & Pricing Management*.

# 1 Introduction

The essentiality of transportation in distribution activities of logistic systems makes the vehicle routing problem (VRP) a combinatorial question of great importance. The VRP originated in 1959 as an integer programming problem (Dantzig, and Ramser, 1959) and aims to find the optimal delivery plan for a fleet of vehicles serving a number of customers. The most common forms of this problem include a central depot and a number of vehicles required to deliver orders to customers at minimum cost. Apart from technical imperatives and operational constraints, the possible sequences of service in VRP increase exponentially with the total number of customers. Thus, when more customers are added to the problem, the computation required to find the solution takes an exponentially longer time. Consequently, acquiring the optimal solution to a VRP is an NP-complete problem (Renaud et al., 1996). Hence, the researchers are interested in developing methods to find high quality solutions from realistic modelling approaches to specific operational constraints of transportation challenges (Lahyani, 2014).

Equally essential to major issues of logistics are the questions of sustainability for transportation operations. Sustainable logistics requires consideration of environmental issues as well as economic efficiency. Research on green transportation is gaining more importance due to the severity of environmental concerns and the undeniable part of transportation activities in them. Transportation is the most substantial factor in depletion of energy resources on Earth. As such, one of the motivations for analyzing VRP stems from the necessity of adopting sustainable practices in transportation planning. To illustrate this point one may consider fossil fuels and their precise utilization. It can make a significant difference not only in environmental imperatives such as control of greenhouse gas (GHG) emissions and global warming, but also in economic performance of transportation systems under new environmental regulations such as carbon taxes. Moreover, the responsible use of energy resources, propagated by environmental campaigners as the first step towards mitigation of air pollution, requires organizations to amend their transportation policies to prevent future catastrophic events. Therefore, in addition to being an optimization problem of high complexity, the VRP is essential from a sustainability viewpoint.

Vehicle routing is also a problem of significant financial importance. Considering economic issues as another aspect of sustainability, VRP deals with a crucial point of industry where a slight improvement in productivity can have far reaching effects on monetary saving from organizational financial resources to national budgets. It is worthy of mention that approximately one tenth of the cost of a finished product is attributed to the costs of transportation activities of its production life cycle (Akerman et al., 2000). Moreover, according to the annual State of Logistics report in 2014, an amount of money equal to 5% of the US gross domestic product is spent on their transportation activities. 77.2% of this amount is attributed to trucking-related activities rather than other modes of transportations (Gilmore, 2014). Such statistics place transportation planning in a crucial niche to be investigated by analytical models of optimization.

Transportation systems' efficiency enhancement is the main objective of green VRP which aims to control pollutant factors and optimize routing. The fuel cost of diesel vehicles makes up a large proportion of total cost in green VRPs (Xiao et al., 2012). As road transportation is the primary source of carbon dioxide emission (Bektaş, and Laporte, 2011), the reduction of emitted GHG is also considered as a part of the objective function in such problems (Erdoğan, and Miller-Hooks, 2012). Transport plays a crucial role in economic development, yet it is the largest consumer of energy resources and the most important factor in global pollution. Therefore, there is an inherent tradeoff within the approaches taken by VRP researchers to prioritizing environmental and economic factors. Eco-friendly policies obtained from green optimization models, in most cases, contradict the optimal solution of classical models comprising of one-sided economic objectives. The green optimization models have been developed under the incentive of aggravating air pollution, to incorporate the GHG emissions produced by different sources in a comprehensive delivery planning. This makes a novelty in approaching VRP with green imperatives. The mobile sources and, particularly, the road transportation are substantial causes of air pollution, overall in respect to carbon dioxide, nitric oxide and volatile organic compound emitted. Therefore, comparative investigation of strategies to reduce the adverse impacts of transportation activities on the environment is a priority for research.

Due to the indispensable role of freight transportation optimization in reducing environmental pollutants (Bauer et al., 2010), it is important to reconsider different parameters of VRP objective functions for gaining effective solutions to the problem. In what follows, different approaches of modifying the VRP explored by contemporary researchers are delineated after a brief background discussion on the origins of such problems.

## 2 Literature Review

An evolution from rudimentary forms of the problem to a variety of sophisticated models associated with different assumptions is evident in VRP literature. Analysis of a large scale Traveling Salesman Problem (TSP) was believed to be the precursor of VRP (Dantzig et al., 1954). However, the first research to tackle a problem with multiple vehicles was investigated in 1964 (Clarke, and Wright, 1964) and the first appearance of the exact terminology goes back to 1977 (Golden et al., 1977). Having the basic concepts institutionalized, Golden introduced probabilistic models of VRP (Golden, and Stewart, 1978) that were extended to uncertain vehicle routing models. Common parameters of uncertainty in VRP were customer demand, travel time, and cost (Gendreau et al., 1996; Lecluyse et al., 2009; Mendoza et al., 2010; Rei et al., 2010). A specific type of the problem, referred to as the Solomon Problems, assumes time windows for serving the customers, incurring a penalty if the product is not delivered within the predefined time window (Solomon, 1987). For a review of algorithms developed to solve VRP with time windows, one may refer to (Bräysy et al., 2004). Another specific type of the problem deals with situations where the information may change through the execution of the delivery. Dynamic VRP models can be reviewed in (Pillac et al., 2013).

It is common for the conceptual formulation of basic VRP models to be built on the distance between customers. However, ignoring time and speed restricts the models from being realistic when observed from a transportation planning viewpoint. In an attempt to remove this issue, Malandraki introduced time-dependent VRP models in which speed varies according to the time of day (Malandraki, and Daskin, 1992). Time-dependent VRP was then investigated more rigorously by (Soler et al., 2009), who incorporated time windows into the model. More realistic ramifications of this type of problem were recently developed by (Hashimoto et al., 2010) (Kritzinger et al., 2012), and (Kok et al., 2012). They focused on hard and soft time windows, impacts of traffic information, and congestion avoidance in time-dependent VRP with time windows. Time and speed were required by these ramifications for measurements such as fuel consumption and emissions of GHG. So, the recent approach contributed to the research area by providing a foundation for investigating not only real world transportation systems, but also green concepts.

Although green VRP has many more aspects to be investigated, only the two mentioned are to be discussed with respect to their relation to the present study. The other aspects, including VRP in reverse logistics, waste collection, end-of-life goods collection, and simultaneous distribution and collection, can be reviewed in (Govindan et al., 2015; Pokharel, and Mutha, 2009; Sbihi, and Eglese, 2007).

Research on vehicle routing with fuel consumption efficiency was quite limited in comparison to classical VRP. Such models were developed to account for speed, load, and distance, as three main factors of fuel consumption, to help obtain effective solutions. Kara et al. were the first to investigate an energy minimizing VRP (Kara et al., 2007). In their suggested model, links were associated with load weight cost in addition to distance cost though there was no formulation provided for fuel consumption. A similar, basic load weight assumption without fuel consumption formulation for a multi-depot problem was later investigated in (Zhang et al., 2011). The idea of formulating fuel consumption was previously suggested by (Sambracos et al., 2004). Later, it was studied by subsequent research articles such as (Marasš, 2008). In addition to load weight cost and distance cost, Kou incorporated speed as an additional factor of fuel consumption into modelling a time-dependent VRP (Kuo, 2010). Xiao suggested a linear function of vehicle load to be embodied in the objective function of a fuel consumption VRP (Xiao, Zhao, Kaku, and Xu, 2012).

Minimizing the emissions of GHG was another green approach to the study of VRP. GHG emissions were implicitly included in the minimization of total distance travelled in classical VRP models. However, responding to transportation challenges encouraged researchers to develop models through which investigation of GHG emissions is more precise. The models that address GHG emissions were built upon the technical building blocks of green transportation, including the pollutant emission estimation method of (Pronello, and André, 2000), the technical report on carbon dioxide emission by (McKinnon, 2007), and the truck freight transportation external costs estimation by (Forkenbrock, 2001). Palmer integrated GHG emission, travel time, and travel distance into a model designed to investigate the impact of vehicle speed on GHG emissions. This model resulted in a potential decrease of 5% in emitted gases (Palmer, 2007). This approach was continued by Sbihi and Eglese, who studied the impacts of traffic on fuel consumption. They focused on the idea that if the engine works at the optimal rotation per minute, the GHG emitted would decrease.

Contradiction between green policies and economic productivity was expected. Improvement in GHG emissions in such models came hand-in-hand with longer routes and slower service (Sbihi, and Eglese, 2007). From a similar approach, Maden varied speed according to the time of day in a duration minimizing model that resulted in a potential decrease in GHG emissions of 7% (Maden et al., 2010). Varying speed in different scenarios was then studied by (Fagerholt et al., 2010), who optimized fuel consumption and emitted GHG in a model with time windows. The first sophisticated research paper to explicitly minimize the GHG emitted was developed by (Ubeda et al., 2011), considering both economic and environmental objectives. It was evident according to the numerical results that using larger freight vehicles can reduce GHG emissions. Multi-objective models of the green VRP was later investigated by (Faulin et al., 2011). Noise pollution was incorporated in their model alongside air pollution and the total distance, the previously suggested terms of the objective function. In a more recent research study, Longo suggested a comprehensive simulation method that optimizes the routes in a sustainable supply chain design model to reduce carbon dioxide as one of the three sustainability aspects taken into account (Longo, 2012).

## 3   Notation and Problem Statement

This section discusses the assumptions and other details of the problem to be investigated in this study. Consider that a fleet of vehicles is going to serve a number of customers in predefined time windows. The start point and the finish point of the routes are both warehouses, and vehicles are limited to load constraints. Different speed limits are assumed with respect to different times of day to incorporate traffic regulations into the problem. In what follows, the notation used for the mathematical formulations are delineated. O.W. is used as shorthand for otherwise. Indices, parameters and variables are as follows.

Sets and indices

| | |
|---|---|
| $V = \{v_0, v_1, \ldots, v_{n+1}\}$ | Set of vertices |
| $A = \{(i,j): i,j \in V, i \neq j\}$ | Set of edges |
| $C = \{v_1, v_2, \ldots, v_n\}$ | Set of customers |
| $K = \{1,2,\ldots,k\}$ | Set of available freight vehicles |
| $\Re = \{1,2,\ldots r\}$ | Set of different speed levels |

Parameters

| | |
|---|---|
| $c_f$ | Constant cost of fuel |
| $e$ | Cost of one gram GHG emitted |
| $\alpha_{ij}$ | Edge constant coefficient |
| $d_{ij}$ | Distance between customer $i$ and $j$ |
| $w$ | Weight of unloaded vehicle |
| $\beta$ | Freight vehicle constant coefficient |
| $[l^r, u^r]$ | Bounds of speed for speed level $r$ |
| $\bar{v}^r$ | Average speed for speed level $r$ |
| $q_i$ | $i$th customers demand |
| $q_{max}$ | Vehicle maximum capacity |
| $[a_i, b_i]$ | Time window for serving $i$th customers |
| $g_i$ | Service duration for $i$th customers |

Decision Variables

| | |
|---|---|
| $x_{ij}^k$ | Equals to one if vehicle $k$ passes the distance between customer $i$ and $j$ O.W. equals to zero |
| $f_{ij}$ | The load carried in edge $i$ to $j$ |
| $y_i^k$ | Start time for serving $i$th customers by vehicle $k$ |
| $z_{ij}^r$ | Equals to one if the freight vehicle passes from $i$ to $j$ with a speed within $[l^r, u^r]$ O.W. equals to zero |

Consider $G = (V, A)$ as a directed graph comprising of a set of vertices $V = \{v_0, v_1, \ldots, v_{n+1}\}$ and a set of edges $A = \{(i,j): i,j \in V, i \neq j\}$. $v_0$ and $v_{n+1}$ represent the warehouses where the freight vehicle of $q_{max}$ capacity is placed. The other vertices are representative of customers. So, the set of customers is a subset of vertices excluding the warehouses $C = \{v_1, v_2, \ldots, v_n\}$. Each customer $v_i$ has its demand $q_i$ and service duration $g_i$ while obviously these two parameters equal zero for each of two warehouses. Customer $i$ expects to receive the service within a specified time window $[a_i, b_i]$. The distance between customer $i$ and customer $j$ are quantified by $d_{ij}$. Finally, upper and lower speed bounds $[l_i, u_i]$ are assumed for a vehicle passing through each edge.

According to the speed and total weight, each freight vehicle emits a certain level of GHG when it passes through an edge. Each gram of GHG emitted is associated with an approximated cost equal to *e* with respect to its environmental consequences. Although GHG emissions are also subjected to parameters such as road slope and gravitational acceleration, they can be controlled by specifying factors such as speed and load carried. This model assumes that all the customers will be served.

## 4 Mathematical Model Formulation

GHG emissions, represented by $E$, are in a direct relationship with the fuel consumption rate, represented by $F$, so a linear function (4.1) is deployed to calculate it according to the GHG-specific emission index parameters $\delta_1, \delta_2$.

$$E = \delta_1 F + \delta_2 \tag{4.1}$$

In contrast to the simplicity of relationship between $E$ and $F$, the fuel consumption rate itself is difficult to calculate. Barth and Boriboonsomsin suggest an approximation formula with eight parameters as shown in (4.2). In the suggested formula, $K$ stands for the engine friction factor, $N$ for engine speed, $V$ for engine displacement, $P_t$ for the tractive power requirement in watts, $\varepsilon$ for the combined efficiency of the transmission and final drive, $P_a$ for engine power demand, $\eta$ for the indicated engine efficiency, and $U$ for a constant coefficient (Barth, and Boriboonsomsin, 2009).

$$\frac{dF}{dt} \approx \left(kNV + \frac{P_{engine}}{\eta}\right) U \ , \quad P_{engine} = \frac{P_t}{\varepsilon} + P_a \tag{4.2}$$

In (4.2) $P_a$ can be expressed as a function of $N$. Similarly $\varepsilon$ can be expressed in terms of $N$ and $P_t$. Therefore, the power requirement on the engine ($P_{engine}$) is a function of the tractive power requirement ($P_t$). Fuel consumption rate is directly related to $P_{engine}$, making it dependent on $P_t$.

The tractive power requirement is dependent on different parameters such as total weight, vehicle speed, and road slope. So, (4.3) is proposed by (Maden, Eglese, and Black, 2010), which uses $v_{ij}$ as the vehicle speed, $w_{ij}$ as the vehicle weight, $f_{ij}$ as the load weight alongside edge constant coefficient $\alpha_{ij}$, and freight vehicle constant coefficient $\beta_{ij}$ to approximate power requirement in edge $(i,j)$.

$$F \approx P_t\left(\frac{d_{ij}}{v_{ij}}\right) \approx \alpha_{ij}(w_{ij} + f_{ij})d_{ij} + \beta v_{ij}^2 d_{ij} \tag{4.3}$$

Therefore, the GHG emission in edge $(i,j)$ can be calculated by $\alpha_{ij}(w_{ij} + f_{ij})d_{ij} + \beta v_{ij}^2 d_{ij}$. The total GHG emission in the network can be incorporated in the objective function as in (4.4).

The latter equation is used as the foundation of modelling GHG emissions in what follows. Let us continue to define the optimization model by introducing the objective function. As demonstrated in (4.4), it aims to minimize the fuel consumption and GHG emission. The capacity of vehicles to be refueled is considered in (4.5). Equation (4.6) guarantees that all of the customers are served in the transportation model. Obviously, the vehicles do their next move from the same customer whose service is just finished as formulated in (4.7). Warehouses are the start point and the finish point of every route and the vehicles arrive to and depart from the warehouses only once. These two common routing principles are mathematically stated in (4.8) to (4.10). Time window constraints are embodied in (4.11) and (4.12), requiring both the start and finish times of the service to be within a predefined range of time.

Another constraint is formulated in (4.13), arguing that the travel time between two nodes has to be within the service time of two customers. This inequality needs more explanation, which will be addressed later in this research. Equation (4.14) balances the network flow and (4.15) guarantees that the load constraint of the vehicle is not violated. Furthermore, equation (4.16) associates a speed to each route. Finally, the types of variables are determined in (4.17) to (4.19).

$$\min \sum_{k} \sum_{(i,j) \in A} (c_f + e)\alpha_{ij} d_{ij} w x_{ij}^k + \sum_{k} \sum_{(i,j) \in A} (c_f + e)\alpha_{ij} f_{ij} d_{ij} \qquad (4.4)$$
$$+ \sum_{k} \sum_{(i,j) \in A} (c_f + e) d_{ij} \beta \left( \sum_{r \in \Re} (\overline{v}^r)^2 z_{ij}^r \right)$$

s.t.

$$\sum_{i \in c} q_i \sum_{j \in v} x_{ij}^k \leq q_{max} \qquad \forall \, k \in K \qquad (4.5)$$

$$\sum_{k} \sum_{j} x_{ij}^k = 1 \qquad \forall \, i \in C \qquad (4.6)$$

$$\sum_{i} x_{il}^k - \sum_{j} x_{lj}^k = 0 \qquad \forall \, l \in C, \forall \, k \in K \qquad (4.7)$$

$$x_{i0}^k = 0, \quad x_{n+1,i}^k = 0 \qquad \forall \, i \in V, \forall \, k \in K \qquad (4.8)$$

$$\sum_{j \in v} x_{0j}^k = 1 \qquad \forall \, k \in K \qquad (4.9)$$

$$\sum_{j} x_{j,n+1}^k = 1 \qquad \forall \, k \in K \qquad (4.10)$$

$$a_i \sum_{j} x_{ij}^k \leq y_i^k \qquad \forall \, i \in V, \forall \, k \in K \qquad (4.11)$$

$$b_i \sum_{j} x_{ij}^k \geq y_i^k \qquad \forall \, i \in V, \forall \, k \in K \qquad (4.12)$$

$$x_{ij}^k \left( y_i^k + g_i + \sum_{r} \left( \frac{d_{ij}}{\overline{v}^r} \right) z_{ij}^r \right) \leq y_j^k \qquad \forall \, (i,j) \in A, \forall \, k \in K \qquad (4.13)$$

$$\sum_{j} f_{ji} - \sum_{j} f_{ij} = q_i \qquad \forall \, i \in C \qquad (4.14)$$

$$q_j x_{ij}^k \leq f_{ij} \leq (q_{max} - q_i) x_{ij}^k \qquad \forall \, (i,j) \in A \qquad (4.15)$$

$$\sum_r z_{ij}^r = \sum_k x_{ij}^k \qquad \forall\, (i,j) \in A \qquad (4.16)$$

$$x_{ij}^k \in \{0,1\} \qquad \forall\, (i,j) \in A \qquad (4.17)$$

$$f_{ij} \geq 0 \qquad \forall\, (i,j) \in A \qquad (4.18)$$

$$z_{ij}^r \in \{0,1\} \qquad \forall\, (i,j) \in A\,, r \in \Re \qquad (4.19)$$

As noted earlier, inequality (4.13) requires a linearization technique. According to the method in (Cordeau et al., 2007) it can be linearized to (4.20) in which $M_{ij}^k$ is calculated with respect to (4.21) to (4.23).

$$y_i^k - y_j^k + g_i + \sum_r \left(\frac{d_{ij}}{\overline{v}^r}\right) z_{ij}^r \leq M_{ij}^k (1 - x_{ij}^k) \qquad \forall\, i \in V\,, j \in C, i \neq j, \forall\, k \in K \qquad (4.20)$$

$$M_{ij}^k = \max\left\{0, b_i + s_i + \frac{d_{ij}}{l} - a_j\right\} \qquad \forall\, i \in V\,, j \in C, i \neq j\,, \forall\, k \in K \qquad (4.21)$$

$$y_j^k + g_i - s_j + \sum_r \left(\frac{d_{j0}}{\overline{v}^r}\right) z_{j0}^r \leq L(1 - x_{j0}^k) \qquad \forall\, i \in V\,, j \in C, i \neq j\,, \forall\, k \in K \qquad (4.22)$$

$$s_j = \left(y_j^k + t_j + \frac{d_{j0}}{v_{j0}}\right) x_{j0}^k \qquad \forall j \in C, \forall\, k \in K \qquad (4.23)$$

The developed mathematical model needs an appropriate solution method to deal with the combinatorial complexity discussed. In the next section, after discussing the methodological background, a solution method is introduced and then tested.

## 5 Solution Method

### 5.1. Methodological background

The noticeable heterogeneity of VRPs and the inherent complexity involved in solving them led to the advent of different tailored algorithms. These algorithms are designed to obtain high quality and near-optimal solutions with computational advantages to exact methods. Exact methods for solving different classes of VRPs including: column generation (Liberatore et al., 2011); branch-and-price (Feillet, 2010); and, a two phase algorithm (Hernandez et al., 2014); are discussed in a recent dissertation (Roberti, 2013). According to (Cordeau, Laporte, Savelsbergh, and Vigo, 2007), non-exact solution methods developed for VRPs include basic heuristics and tailored meta-heuristics. Basic heuristics are comprised of: a saving heuristic; a sweeping heuristic; and, a Fisher-Jaikumar heuristic. Tailored meta-heuristics developed consist of: genetic algorithms; tabu searches; and simulated annealing algorithms. The ability of such algorithms to acquire a high quality solution, referred to as accuracy, needs to be individuated quantitatively. Besides being accurate, a well-performing evolutionary algorithm constructs a basic, good quality solution and defines appropriate neighborhood selection strategies to enhance the solutions.

The solution method used in this study was first developed as a calculation method for computing machines by (Metropolis et al., 1953), which then evolved into Simulated Annealing (SA) by (Kirkpatrick et al., 1983) as a meta-heuristic method for combinatorial optimization. According to the algorithm, the search procedure starts by expanding basic solutions into new solutions. The objective function values are conceptualized as energy levels. From this point, a controlling parameter, known as the temperature, is decreased according to a predefined function. The algorithm continues by producing a new population and assessing it. In each of the iterations, the previous solutions are replaced by new solutions with a larger fraction of high quality ones. Continuing the search beyond local optimality is guaranteed by considering the probability of replacing a solution with another solution of lower quality. Similarly, a gradual decrease of temperature ensures the convergence of the algorithm to an appropriate solution when the controlling parameters are meticulously tuned.

The algorithm in this study includes four different neighborhood selection strategies to expand the new solutions within different directions of the feasible space. The temperature is reduced when the number of moves reaches a certain point. In contrast to the other meta-heuristic methods, only enough neighborhoods are produced to obtain one new feasible solution. The flowchart in Fig. 1 demonstrates how the suggested SA works.

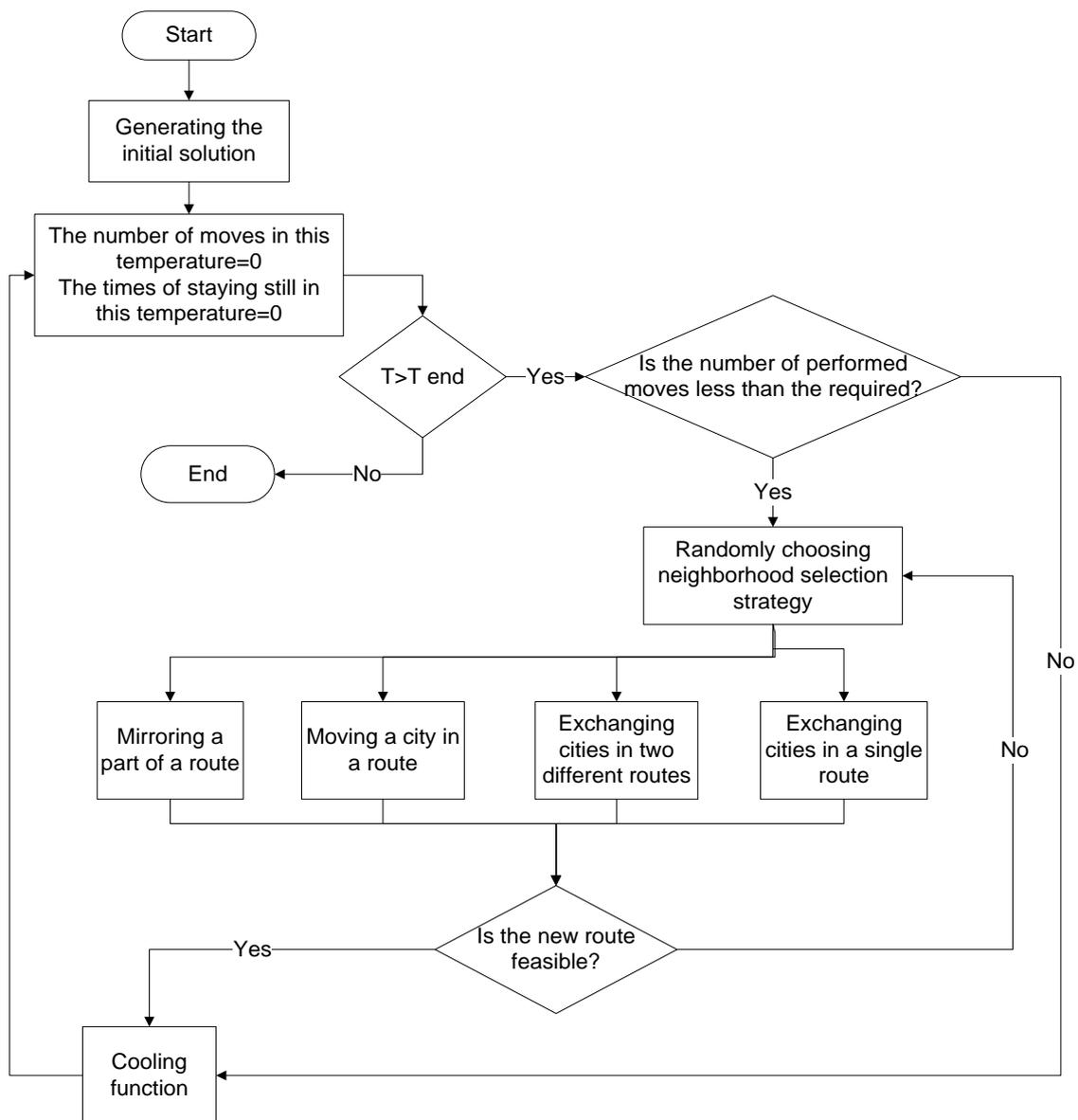

**Fig. 1 Simulated annealing flowchart drawn by Microsoft Visio**

## 5.2. Representation and neighborhood selection

As the representation of the solutions affects the computational complexity, it is important to represent the solutions as simply as possible. The routes and their associated speed limits are the two pieces of information incorporated in a string to characterize a solution. The string 0,1-3,1-5,1-7,1-8,1-11,1-0,1-2,1-4,1-6,1-11,1-0,2-1,2-9,2-10,2-11,2 represents a solution comprising of three routes. Each route is planned to allow for a vehicle to drive within the range of the predefined average speeds according to the time of day. A comma mark is used to distinguish the cities from the average speeds. The zeros and elevens in the string represent the warehouses at both start and finish points of each route. The numbers after a comma and before a hyphen stand for the average speed limits. So in the exemplified string, the first two vehicles are scheduled to drive within the first average speed bracket and the last one is planned to reach the customers by driving within the second average speed bracket. The remaining numbers in the string represent the order of the customers to visit. The delivery plan of this string is illustrated in Fig. 2.

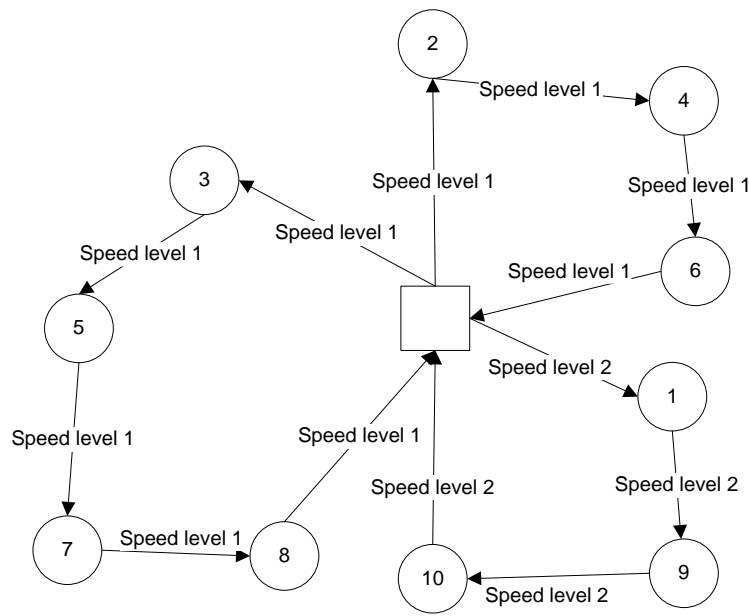

**Fig. 2 Delivery plan for the string exemplified**

An initial node heuristic (Mosheiov, 1994) is employed in this study to establish an initial solution. According to this method, a random permutation of all the customers is first created by the algorithm. Then the time of reaching customer $i$ is checked to ensure it is consistent with the corresponding time window. Within such a procedure, if the two are not consistent, the customer in question is moved to the end of the route and the next customer is moved to the beginning of the route. As previously discussed, four different neighborhood selection strategies are deployed in this study to obtain new solutions by making alternations to current solutions. The first strategy randomly chooses a portion of the solution and replaces it with its mirror. The second strategy deals with obtaining a new solution by moving a customer within the sequence. The exchange of two customers in a solution is the foundation of the third strategy. Finally, the fourth strategy takes a portion of the route from a solution and replaces it in a randomly selected place in another solution. These four functions help to move within the feasible space and continue the search. Fig. 3 illustrates how the four neighborhood selection functions affect the string exemplified in Fig.2.

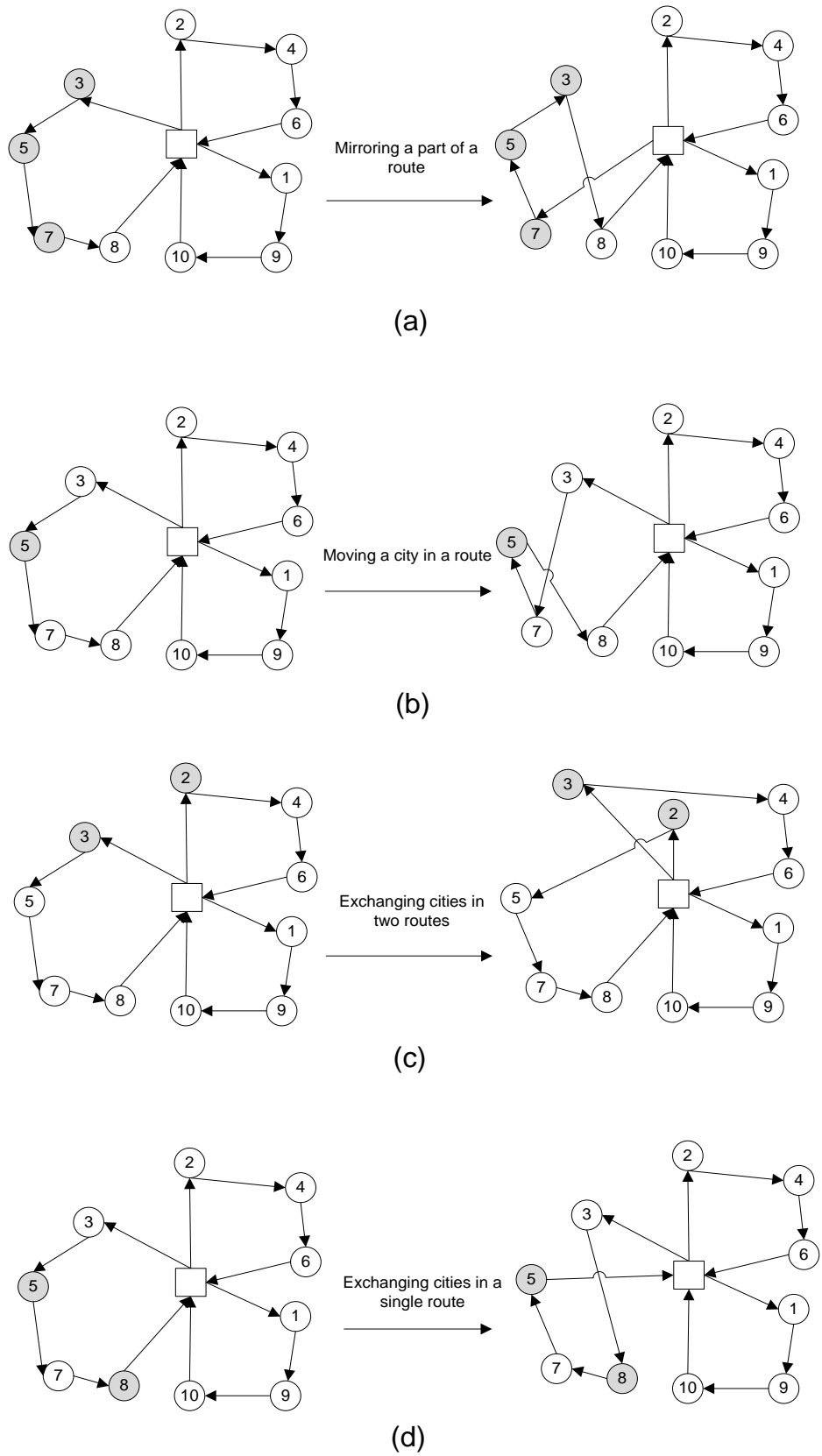

**Fig. 3** Neighborhood selection strategies (a. Mirroring a part of a route, b. Moving a city in a route, c. Exchanging cities in two route, d. Exchanging cities in a single route)

### 5.3. SA Parameters

The meta-heuristic method deployed to solve the problem includes four controlling parameters. Although it is assumed that the algorithm stops when it reaches the final temperature, considering a stopping criterion is required for each temperature. The number of generated neighborhoods is restricted in each iteration. So, the iteration is skipped unless a feasible solution is generated before reaching $n + 50$ in which $n$ represents the number of customers. Moreover, an exponential cooling function with a parameter equal to 0.97 is found by trial and error to account for changing the algorithm stage. Furthermore, the initial temperature is set to 1, as in such a setting at least half of the moves result in a decrease of objective function value for the new solutions developed. It is also made evident that the algorithm does not progress any further after reaching a temperature equal to 0.001, so this value is set as the final temperature.

## 6  Computations and results interpretation

In this section a number of randomly generated test problems of different sizes are solved, ranging from small problems with five customers to large problems comprising of a hundred customers. The empty vehicles are assumed to weigh ten units and be capable of transporting a load equal to their own weight. Moreover, two different time ranges are considered, comprising of a day time range limited to a speed of 60 km/h and another time range for afternoon and night limited to a speed of 50 km/h.

The exact solution method was coded by OPL 6.3 in Cplex 12.2 and the SA algorithm was programmed by MATLAB 2011A. Both programs were run by an Intel Core-i5 computer with 4 GBs of RAM. The results of the computations are summarized in Table 1, including the values of objective functions and computation time for each test problem. The computation times for the exact solution of numerical examples with more than 20 customers were restricted to 30 minutes. Therefore, the respective cells in Table 1 show the best objective function value found in the restricted time.

**Table 1 Computational result for randomly generated test problems**

| No. | Number of Customers | Exact Method | | Meta-heuristic Method | | Gap (%) | Time decreased (%) |
|---|---|---|---|---|---|---|---|
| | | Optimal Value | Computation Time (s) | Best Solution Found | Computation Time (s) | | |
| 1 | 5 | 17551017.83 | 10.43 | 17551017.83 | 10.65 | 0 | -2.1 |
| 2 | 6 | 11080856.63 | 16.76 | 11080856.63 | 9.40 | 0 | 43.9 |
| 3 | 7 | 12503316.03 | 9.25 | 12503316.03 | 9.62 | 0 | -4.0 |
| 4 | 8 | 12515151.04 | 32.98 | 12948361.64 | 9.88 | 3 | 70.0 |
| 5 | 9 | 13034308.25 | 26.56 | 13513451.74 | 10.03 | 3.5 | 62.2 |
| 6 | 10 | 18551017.83 | 28.26 | 18551017.83 | 10.43 | 0 | 63.1 |
| 7 | 11 | 16639092.82 | 56.42 | 16639092.82 | 10.77 | 0 | 80.9 |
| 8 | 12 | 18904523.19 | 245.85 | 19456017.21 | 10.91 | 2.8 | 95.6 |
| 9 | 20 | 25421646.51 | 428.87 | 25421646.51 | 13.12 | 0 | 96.9 |
| 10 | 50 | 61658229.14* | 1800.00 | 62407114.52 | 67.45 | 1.2 | 96.3 |
| 11 | 70 | 102155572.15* | 1800.00 | 104774945.8 | 110.71 | 2.5 | 93.8 |
| 12 | 80 | 121125965.29* | 1800.00 | 119335926.4 | 176.28 | - | 90.2 |
| 13 | 90 | 150161797.56* | 1800.00 | 146499314.7 | 375.28 | - | 79.2 |
| 14 | 100 | 116314816.95* | 1800.00 | 113146709.1 | 165.77 | - | 90.8 |

*Best objective function value obtained by Cplex after 30 min.

As evident in Table 1, except for test problem number 4, 5, and 8, with an insignificant gap between the results of two methods, the SA algorithm was capable of finding the exact solutions for the rest of small to mid-sized test problems. According to the results obtained for large test problems with 50 to 100 customers, the SA algorithm performed well. It even came up with solutions of higher qualities in comparison to the exact method in much shorter computation times.

The algorithm convergence for a test problem with 100 customers is plotted in Fig. 4. As it is shown, the SA algorithm initially resisted the changes due to the quality of the initial solution. Then in the middle of the computations, the solutions evolved into the final value with an almost uniform trend of improvement.

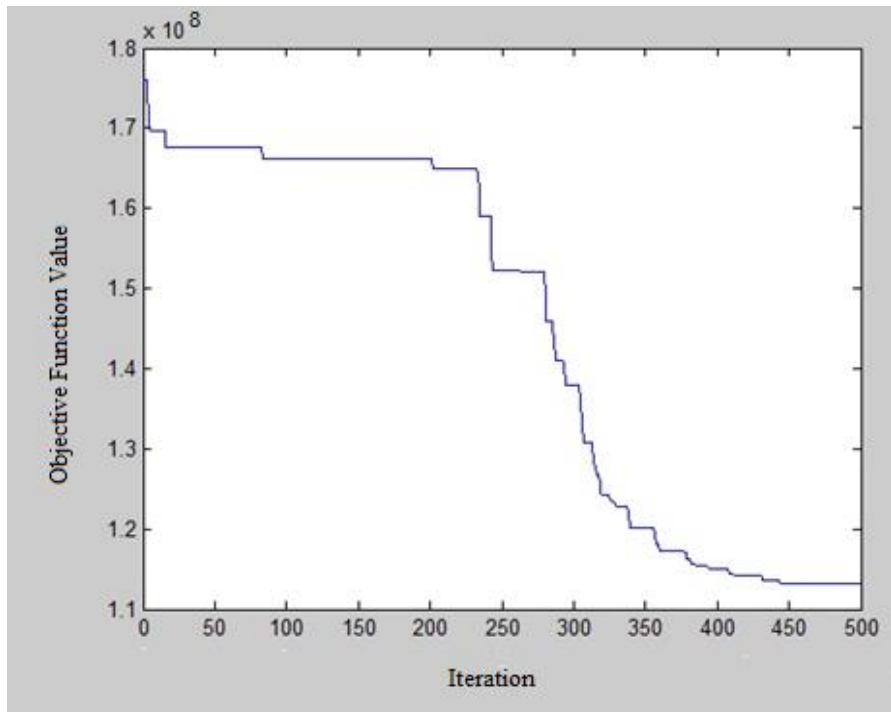

**Fig. 4 The trend of convergence in reaching the final solution for test problem 14**

In order to evaluate the time performance of the suggested algorithm, the computation times of the first nine test problems are compared in Fig. 5.

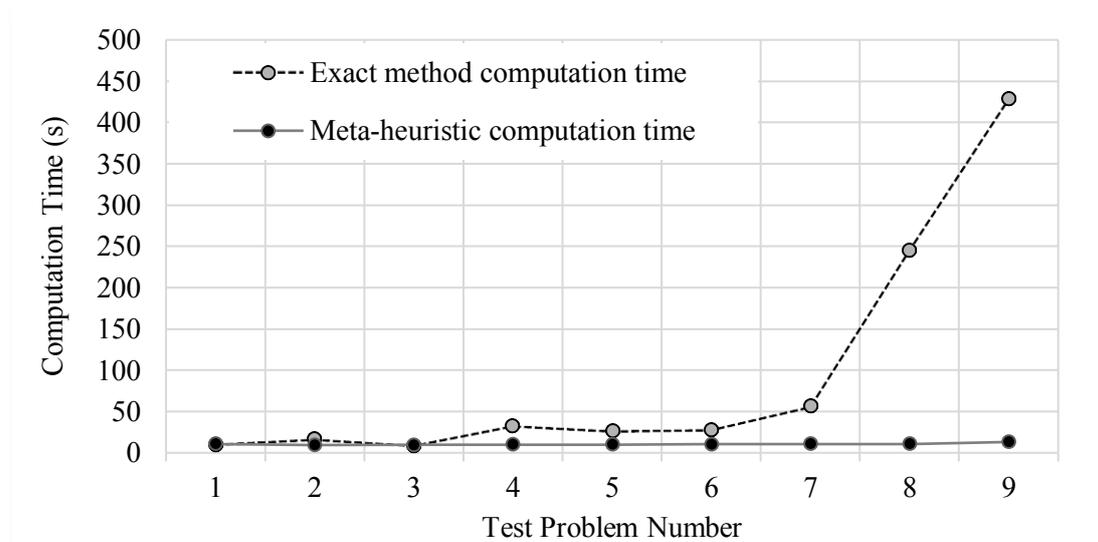

**Fig.5 Computation time performance comparison**

The simulated annealing method almost took a time equal to the exact method computation time for the first three test problems, as evident in Fig. 5. However, the computation times for the next four test problems were significantly reduced. The last two test problems were associated with a prolonged computation time in the exact method that is also significantly reduced to a more appropriate solve time by the suggested solution method.

## 7 Further Discussion

The proposed model is developed based on incorporating the adverse impacts of transportation into the objective function of a minimization model. One may argue that in many cases of transportation planning, classic decisions are made based on the analysis of total cost for different alternatives while environmental imperatives are a secondary consideration. In such cases, meeting the environmental requirements by following normative guidelines has been the standard practice of supply chain managers. In spite of this, the new environmental economics solutions advocated to mitigate environmental problems require new decision making paradigms. For instance, based on environmental taxations introduced in developed countries, optimal level of emissions should allow pollution occurring only if the benefits to society exceeds the costs of pollution (Claus et al., 2010). Our proposed model is capable of providing a tool for this novel decision making paradigm in which total cost is not the only determinant factor.

## 8 Application Perspectives

Although this study focuses on GHG emissions and fuel consumption in VRP, the findings may well have a bearing on sustainable supply chain and logistics as a broader applications. This research has several practical applications. Firstly, the proposed model can be used to extend common location technologies into delivery planning applications that provide sustainable solutions for delivery companies using location data and traffic maps. Secondly, as discussed earlier a slight improvement in delivery productivity can have far reaching effects on monetary savings for courier companies. The proposed model can be deployed as the backbone of a decision making tool for green courier companies to help them save resources and fuel required to deliver the same number of packages per day while environmental aspects of operations are taken into account as well.

## 9 Conclusion and Future Research

The environmental challenges of freight transportation call for numerical models to make sustainably responsible decisions in delivery plans. Fuel consumption and GHG emissions are highly regarded challenges of transportation networks due to their monetary cost and undeniable role in environmental degradation. This study has posed a new problem with more relevance to current environmental concerns. The problem was solved using a combinatorial optimization approach capable of handling industry-sized models.

To be more specific, this study aimed to propose a model for a capacitated VRP with time windows and time-dependent speed limits. Fuel consumption and GHG emissions were minimized as two principal green objectives to investigate the feasibility of planning transportation activities with respect to green imperatives. Accordingly, the control of speed and time of travel were studied through a mixed-integer programming model with constant values for nodes' demand.

Solving the proposed model to optimality using exact methods was an NP-complete problem. So, an SA algorithm capable of solving problems of different sizes was suggested. The proposed algorithm was characterized by simple solution representation and productive neighborhood selection strategies. The SA algorithm was tested in a number of numerical examples. It obtained high quality solutions with efficiency in computations by establishing initial solutions and evolutionarily improving them using different neighborhood selection strategies.

As this study posed a new problem in green transportation, it offered significant research gaps for further investigation. The developed green objective function is capable of being investigated more pervasively in different types of VRPs, as discussed earlier. The proposed model can be extended into a multi-objective optimization model to investigate the trade-offs between the two green objectives and the total cost that can be potentially in conflict. Similarly, one may introduce a model to minimize total cost while controlling the fuel consumption and GHG emissions at certain levels inferred from guidelines. Using the same research structure of this study, other environmental aspects of freight transportation can be considered and formulated as optimization models to explore other research avenues in mitigating impacts of transportation by conducting quantitative analyses. This study has illustrated the potential of planning transportation according to environmental objectives. Therefore, it is necessary to investigate these concepts further in order to advance toward practical implementations.

**Acknowledgements**

*A Green Perspective on Capacitated Time-dependent Vehicle Routing Problem*